\documentclass{article}
\usepackage[margin=.98in]{geometry}
\usepackage{amssymb}
\usepackage{amsmath}
\usepackage[all]{xy}
\usepackage{hyperref}
\usepackage{color,bm}

\numberwithin{equation}{subsection}
\newcounter{defi}
\renewcommand{\thedefi}{\arabic{subsection}.\arabic{defi}}
\numberwithin{defi}{subsection}
\newcommand{\newdefi}[1]{\refstepcounter{defi}   \textbf{#1~\thedefi.}} 
\newcommand{\bewende}{\hfill $\square$\vspace{2mm}\\}
\newcommand{\real}{\mathrm{Re}\,}
\newcommand{\imag}{\mathrm{Im}\,}
\newcommand{\Hcurl}{H_{\mathrm{curl}}}

\makeatletter
\def\blfootnote{\gdef\@thefnmark{}\@footnotetext}
\makeatother

\begin{document}

\title{An inverse problem for Maxwell's equations \\with Lipschitz parameters\blfootnote{Date: \today}}
\author{Monika Pichler}
\date{}
\maketitle

\begin{abstract}
We consider an inverse boundary value problem for Maxwell's equations, which aims to recover the electromagnetic material properties of a body from measurements on the boundary. We show that a Lipschitz continuous conductivity, electric permittivity, and magnetic permeability are uniquely determined by knowledge of all tangential electric and magnetic fields on the boundary of the body at a fixed frequency.
\end{abstract}

\subsection{Introduction}
We consider the time-harmonic Maxwell's equations for the electric field $E$ and magnetic field $H$,
\begin{equation}
 \nabla \wedge E - i\, \omega\, \mu\, H =0, \qquad \qquad \nabla \wedge H + i\, \omega\, \gamma\, E =0, \label{eq:ME}
\end{equation}
in a bounded domain $\Omega \subset \mathbb{R}^3$ with Lipschitz boundary, 
where the constant $\omega > 0$ is the angular frequency; $\gamma = \varepsilon + i\sigma/\omega$; and the electric permittivity $\varepsilon$, the magnetic permeability $\mu$ and the conductivity $\sigma$ are Lipschitz functions. 
It is known (see e.g. \cite{C2010}) that for a given set of Lipschitz continuous parameters $\mu,\,\varepsilon,$ and $\sigma$, and tangential boundary data in a suitable space ($\nu$ denotes the outward unit normal vector on $\partial \Omega$),
\[
\nu \wedge E = f\in  TH(\partial \Omega) = \left\{ F \in B^{-1/2}(\partial \Omega)^3: \exists u \in \Hcurl(\Omega), \nu\wedge u = F\right \} ~~~\mathrm{on}~\partial \Omega, 
\]
the equations \eqref{eq:ME} have a unique solution pair $(E,H)$ in the space $\Hcurl(\Omega)\times \Hcurl(\Omega)$, except for a discrete set of values of $\omega$. Here,
\[
\Hcurl(\Omega) = \big\{ u\in L^2(\Omega)^3: \nabla \wedge u \in L^2(\Omega)^3 \big\},
\]
and $B^{-1/2}(\partial \Omega)$ is a Besov space on $\partial \Omega$, see \cite{C2010}.

We study the inverse problem, which aims at recovering the parameters $\mu,\,\varepsilon$, and $\sigma$, given some information about solutions on the boundary of the domain. More precisely, for a given frequency $\omega$, suppose we know the \emph{Cauchy data set} of tangential boundary values of the electric and magnetic fields $E$ and $H$,
\[
C(\mu,\varepsilon,\sigma;\omega)=\big\{ (\nu\wedge E|_{\partial \Omega},\nu\wedge H|_{\partial \Omega}): (E,H)~\mathrm{solves~\eqref{eq:ME}~with~parameters~}\mu,\varepsilon,\sigma \big\}.
\]
The question we want to answer is whether this set of data uniquely determines $\mu,\,\varepsilon$, and $\sigma$ in $\Omega$. We are going to prove the following result.\vspace{2mm}

\newdefi{Theorem} \label{thm1} \emph{Let $\Omega \subset \mathbb{R}^3$ be a non-empty bounded domain such that $\partial \Omega$ is locally described by the graph of a Lipschitz function. Fix $\omega>0$, and let $\mu_1,\,\varepsilon_1,\,\sigma_1,\, \mu_2,\,\varepsilon_2,\,\sigma_2\in C^{0,1}(\overline{\Omega})$ be bounded Lipschitz functions such that for a positive constant $c_o$ and $j=1,2$,
\[ 0<\mu_o \leq \mu_j(x),~~~ 0< \varepsilon_o \leq \varepsilon_j(x),~~~ 0\leq \sigma_j(x) ~~~\forall x\in \overline{\Omega},\]
\[|\mu_j(x)-\mu_j(y)| \leq c_o|x-y|,~~~|\varepsilon_j(x)-\varepsilon_j(y)| \leq c_o|x-y|,~~~|\sigma_j(x)-\sigma_j(y)| \leq c_o|x-y|~~~\forall x,y \in \overline{\Omega}. \]
Assume further that $\mu_1(x)=\mu_2(x)$, $\varepsilon_1(x)=\varepsilon_2(x)$, and $\sigma_1(x)=\sigma_2(x)$ for all $x\in \partial \Omega$. Then} 
\[ C(\mu_1,\varepsilon_1,\sigma_1;\omega)=C(\mu_2,\varepsilon_2,\sigma_2;\omega)~~ \Rightarrow ~~\mu_1=\mu_2,
 ~~\varepsilon_1=\varepsilon_2~~{and}~~\sigma_1=\sigma_2. \]

Questions of this kind have been extensively studied for a number of different equations. In 1980, A.P. Calder\'{o}n \cite{C1980} posed the question whether one could recover the conductivity of a body from measurements of the electric voltage and current on the surface of the body. This is the inverse problem for the conductivity equation,
\[ \nabla \cdot (\sigma \nabla u)=0~~\mathrm{in}~\Omega,\]
where $\sigma$ is the conductivity and $u$ is the electric potential in a bounded body $\Omega \subset \mathbb{R}^n$. The goal is to recover $\sigma$ from the knowledge of the {\it Dirichlet-to-Neumann map} (DN map) $\Lambda_\sigma$, which formally maps $u|_{\partial \Omega}$ to $ \nu\cdot \sigma \nabla u|_{\partial \Omega}$, and is defined by duality as
\[ \big\langle \Lambda_\sigma f,g\big\rangle = \int_{\Omega} \sigma \nabla u \cdot \nabla v dx, \]
where $u$ solves the conductivity equation with $u|_{\partial \Omega}=f$, and $v|_{\partial \Omega}=g$.

For scalar conductivity $\sigma$, this question was answered affirmatively in \cite{SU1987} for a smooth conductivity on a smooth domain in dimension greater than or equal to three. The authors assumed two conductivities had identical DN maps and derived an integral identity involving the difference of the conductivities as well as a product of solutions to the respective conductivity equations; they then constructed a special family of solutions to the conductivity equation, so-called {\it complex geometrical optics (CGO) solutions}, that are of the form
\[ w(x)= e^{\zeta \cdot x} \big( 1 + r_{\zeta}(x) \big), \]
with $\zeta \in \mathbb{C}^n$ such that $\zeta\cdot \zeta=0$, and $r_{\zeta}(x)$ becoming small in a suitable sense as $|\zeta|$ becomes large. These solutions, when plugged into the integral identity, allowed to conclude that the two conductivities in fact had to be equal. Since this seminal work, the method developed there has been modified in various ways, allowing to lower the required smoothness of the domain and conductivity in several instances. A survey of results concerning the conductivity equation can be found in \cite{U2009}. 
We want to highlight two papers in this context which we will repeatedly refer to in what follows: in \cite{HT2013}, uniqueness for continuously differentiable conductivities on a Lipschitz domain was shown, using an estimate holding \emph{on average} that allowed to do away with the low regularity of the conductivity, and an adaptation of this estimate will be used also in this paper. In \cite{CR2016}, uniqueness has been shown for Lipschitz conductivities on a Lipschitz domain. We will also heavily rely on a modification of a Carleman type estimate from this paper to show existence of solutions to some auxiliary equations.

Estimates of Carleman type, involving different types of weight functions, have been employed extensively in the study of partial data problems, in which the DN map is known only on part of the boundary, to construct CGO solutions with controlled behavior on parts of the boundary. In \cite{BU2002}, Carleman estimates with a linear weight function were used to construct CGO solutions and show unique determination of a twice continuously differentiable conductivity, if the DN map is known on slightly more than half of the boundary. In \cite{KSU2007}, nonlinear weights are used to show that the DN map measured on an arbitrary open subset of the boundary uniquely determines a twice continuously differentiable conductivity. A constructive proof of this result which required choosing particular uniquely specified CGO solutions was given in \cite{NS2010}. A different treatment of partial data problems was presented in \cite{I2007}, where a reflection argument was used, assuming the inaccessible part of the boundary was part of a sphere or plane.

The above results concern spatial dimensions greater than or equal to three. In the two-dimensional case, the problem needs to be treated differently, and tools from complex analysis are employed. In this setting, the first global uniqueness result is \cite{N1996} for in a suitable sense twice differentiable conductivity; the regularity requirement in this case could be lowered significantly to bounded measurable conductivities in \cite{AP2006}.

We also want to point to some boundary identifiability results. For the conductivity equation, it was shown in \cite{KV1984} that knowledge of the DN map fully characterizes a smooth conductivity on the boundary; the corresponding result for Lipschitz conductivities on a domain with Lipschitz boundary was proved in \cite{A1990}. This is particularly important, since the inverse problem for the conductivity equation is typically treated by transforming it to a related problem for a Schr\"{o}dinger equation, and knowledge of the conductivity and its normal derivative on the boundary is necessary for this process. A similar result for smooth parameters exists for Maxwell's equations: it was shown in \cite{JM2000} and \cite{M1997} that knowledge of the Cauchy data set fully determines smooth electromagnetic parameters on the boundary of a smooth domain. 

The method of using CGO solutions introduced in \cite{SU1987} has been adapted to facilitate studying different equations, see e.g. \cite{NSU1995} for a Schr\"{o}dinger equation with a magnetic potential and \cite{NU1994,NU2003} for an inverse problem in elasticity. 
The inverse problem for Maxwell's equations was first formulated in \cite{SIC1992}. The first global uniqueness result for smooth parameters on a smooth domain is \cite{OPS1993}; another proof was given in \cite{OS1996}. More recently, in \cite{CZ2014} the case of continuously differentiable parameters on a domain with $C^1$ boundary was examined. Partial boundary data problems were studied in \cite{COS2009} using the reflection argument introduced in \cite{I2007}, and in \cite{COST2016}, extending the ideas from \cite{BU2002,KSU2007} to Maxwell's equations. The inverse problem on a manifold has been studied in \cite{OPS2003}, and \cite{KSU2011} considered the problem in a non-Euclidean setting.

The paper \cite{OS1996} introduced the idea of relating Maxwell's equations to a matrix Schr\"{o}dinger equation and using solutions to this elliptic equation to obtain solutions to Maxwell's equations. This approach has been used in many publications since that dealt with Maxwell's equations, and it is also the starting point for our analysis. 
We will modify Maxwell's equations in order to obtain an auxiliary elliptic first order matrix equation, as well as a closely related matrix Schr\"{o}dinger equation. 
The original method in \cite{OS1996} requires at least twice differentiable parameters, so at the present level of regularity, we obtain a Schr\"{o}dinger equation with a weakly defined potential, and all the equations that follow are to be understood in a weak sense. 
We will then derive an integral formula involving special solutions to the auxiliary equations. This formula will be the starting point of the uniqueness proof: we will construct CGO solutions to plug into the integral formula, which are of the form
\[ w(x)= e^{\zeta\cdot x}(A_{\zeta}+R_\zeta(x)), \]
with $A_\zeta$ a consant vector and $R_\zeta$ a vector function, depending on a large parameter $|\zeta|$. We then take the limit of the large parameter in the integral formula, which will yield a set of differential equations for the unknown parameters. In order to be able to perform this limit process, the involved functions need decay properties in suitable norms. If the parameters are at least twice continuously differentiable, one obtains a suitable estimate in weighted $L^2$ spaces (c.f. \cite{SU1987}) for the function $R_\zeta$. In order to deal with the lower regularity in the present situation, we will adapt some tools developed in the context of the conductivity equation: we will use the Carleman estimate from \cite{CR2016}, which allows us to construct solutions the the auxiliary equations with Lipschitz electromagnetic parameters, as well as an estimate derived in \cite{HT2013}, which concerns decay of $R_\zeta$ in an averaged sense in suitable weighted spaces that are adapted to the equations we are working with. It turns out that this decay on average is sufficient to be able to perform the limit in the integral formula. 

The difficulty in our situation when compared to the inverse problem for the conductivity equation is to make sure the CGO solutions to the matrix Schr\"{o}dinger equation provide solutions to Maxwell's equations; this requires uniqueness of solutions in some sense that guarantees we obtain suitable CGO solutions. The construction of CGO solutions using a Carleman estimate however does not in general provide uniqueness. In order to insure uniqueness, we were inspired by the approach used in \cite{NS2010}, where uniquely determined solutions were needed to be able to reconstruct the conductivity. 
More precisely, a Carleman estimate allows to construct a bounded functional on a subspace of the solution space, which in \cite{NS2010} is $L^2(\Omega)$, while in our case, we are dealing with a Bourgain-type space. We choose the space carefully to guarantee uniqueness of the solution.

This paper is organized as follows: In Section \ref{sec:aux}, we will derive the auxiliary elliptic equations we will work with; we also present some results on how solutions to the Schr\"{o}dinger equation are related to those of the auxiliary first order equation. In Section \ref{sec:intform}, we obtain the integral formula involving the unknown parameters and solutions to the auxiliary equations.  
We will adapt the a priori estimate from \cite{CR2016} to our situation in Section \ref{sec:apriori} and use it in Section \ref{sec:cgo} to construct the CGO solutions, and go on to show they satisfy an averaged estimate in the spirit of \cite{HT2013}. Finally, we will perform the uniqueness proof in Section \ref{sec:unique}, by deriving differential equations for the parameters as outlined above, and showing a uniqueness result for these equations, from which unique solvability of the inverse problem follows.
\vspace{2mm}

\subsection{An auxiliary elliptic system}\label{sec:aux}

As in \cite{OS1996}, we start by augmenting Maxwell's equations to obtain an elliptic system. Assume $E$ and $H$ satisfy \eqref{eq:ME}, and 
define the scalar potentials (which in this case vanish identically)
\begin{equation}
\Phi = \frac{i}{\omega} \nabla \cdot \left( \gamma E \right), \qquad
\Psi = \frac{i}{\omega} \nabla \cdot \left( \mu H \right). \label{eq:phipsi} 
\end{equation}
We use these to modify Maxwell's equations to obtain the equations 
\begin{equation}
\nabla \wedge E - \frac{1}{\gamma} \nabla \frac{1}{\mu} \Psi - i\, \omega\, \mu\, H =0,\qquad 
\nabla \wedge H + \frac{1}{\mu} \nabla \frac{1}{\gamma} \Phi + i\, \omega\, \gamma\, E = 0, \label{eq:MEaug} 
\end{equation}
which at the given level of regularity of the parameters are satisfied in a weak sense.
If we rescale the fields and potentials as
\[ e=\gamma^{1/2} E, \quad h = \mu^{1/2} H, \quad \phi = \frac{1}{\gamma\mu^{1/2}} \Phi, \quad \psi = \frac{1}{\gamma^{1/2} \mu} \Psi, \]
then from equations \eqref{eq:phipsi}-\eqref{eq:MEaug} it follows that the vector $X=(\phi,e,h,\psi)^T$ is a weak solution to the matrix differential equation
\[ \mathcal{P} X := (P(i\nabla) - k + V)X =0, \]
where $P(i\nabla)$ is the elliptic first order matrix differential operator 
\[ P(i\nabla) = i \left(\begin{array}{cccc}
0 &\nabla\cdot &0&0\\
\nabla & 0 &\nabla \wedge & 0 \\
0 & - \nabla \wedge & 0 & \nabla \\
0 & 0 & \nabla\cdot & 0 \end{array} \right),
 \] 
and 
\[ V = (k-\kappa)I_8 + \big(P(i\nabla) D\big) D^{-1},\] 
with $D=\mathrm{diag}(\mu^{1/2}, \gamma^{1/2}I_3,\mu^{1/2}I_3,\gamma^{1/2})$, $k=\omega(\epsilon_o\mu_o)^{1/2}$ and $\kappa=\omega(\mu\gamma)^{1/2}$. We also consider the operator $\mathcal{P}':=
P(i\nabla) +k-V^T$. 
It is important for our analysis to note that 
\[  \mathcal{P}\mathcal{P}' = (P(i\nabla) - k + V)(P(i\nabla) + k -V^T) = -(\Delta +k^2) +Q, \]
where $Q$ is the weakly defined matrix multiplier
\[ Q= VP(i\nabla) - P(i\nabla)V^T+k(V+V^T)-VV^T. \]
Note that the modification has been made in such a way that this operator is a zeroth order operator. With $\alpha=\nabla \log \gamma$, $\beta= \nabla \log \mu$, and $\theta =\omega^2(\gamma\mu-\varepsilon_o\mu_o)$
, $Q$ is of the following form, where we write 8-vectors $w=(w_1,w_2,w_3,w_4)^T$ with scalar functions $w_1,w_4$, and 3-vectors $w_2, w_3$, and $\varphi$ is a vector test function,
\begin{eqnarray}
\big \langle Q w\,,\, \varphi \big \rangle &=&\int \bigg( \frac{1}{4}\big(|\alpha|^2-4\theta \big) \big[ w_1\varphi_1 + w_3\cdot \varphi_3 \big] -\frac{1}{2} \alpha \cdot \Big[ \nabla(w_1 \varphi_1-w_3 \cdot \varphi_3) +\nabla\cdot (w_3 \varphi_3^T + \varphi_3 w_3^T) \Big]\notag \\
&& \qquad + \frac{1}{4}\big(|\beta|^2-4\theta \big) \big[ w_4\varphi_4 + w_2\cdot \varphi_2 \big] -\frac{1}{2} \beta \cdot \Big[ \nabla(w_4 \varphi_4-w_2 \cdot \varphi_2)+\nabla\cdot (w_2 \varphi_2^T + \varphi_2 w_2^T) \Big] \notag\\
&&\qquad - 2 i \kappa \nabla \cdot  \Big(w_1\varphi_2 + w_2 \varphi_1+w_3\varphi_4 + w_4 \varphi_3
\Big)\bigg)dx.\label{eq:Qweak}
\end{eqnarray}
We will also be using the operator
\[  \mathcal{P}'\mathcal{P} = (P(i\nabla) + k -V^T)(P(i\nabla) - k + V) =-(\Delta +k^2) +\tilde{Q}, \]
where $\tilde{Q}$ is defined by
\begin{eqnarray}
\big \langle \tilde{Q} \,w\,,\, \varphi \big\rangle &=& \int \bigg( \frac{1}{4}\big(|\beta|^2-4\theta \big)\big[w_1 \varphi_1 + w_3\cdot \varphi_3\big] + \frac{1}{2} \beta \cdot \Big[ \nabla(w_1 \varphi_1-w_3 \cdot \varphi_3) + \nabla\cdot (w_3 \varphi_3^T + \varphi_3 w_3^T) \Big] \notag \\
&&\qquad + \frac{1}{4}\big(|\alpha|^2-4\theta \big)\big[w_4 \varphi_4 + w_2\cdot \varphi_2\big] + \frac{1}{2} \alpha \cdot \Big[ \nabla(w_4 \varphi_4-w_2 \cdot \varphi_2) +\nabla\cdot (w_2 \varphi_2^T + \varphi_2 w_2^T) \Big] \notag \\
 &&\qquad + 2 i \kappa\nabla \cdot  \big(w_3\wedge \varphi_2- w_2\wedge \varphi_3\big) \bigg) dx.\label{eq:Qtilde}
\end{eqnarray}
It is noteworthy that the first and last components of this operator decouple from the rest and allow us to treat those separately. This will be important later on, as we proceed as follows. We will construct solutions to Maxwell's equations by first finding CGO solutions to the Schr\"{o}dinger equations, and then applying the auxiliary operators $\mathcal{P}$ and $\mathcal{P}'$. If the resulting functions have vanishing first and last components, they give solutions to Maxwell's equations. In showing that these components vanish, the shape of $\tilde{Q}$ will play an important role.

The following propositions establish how solutions to Schr\"{o}dinger equations help us in finding solutions to the auxiliary equations
\[
\mathcal{P} v = (P(i\nabla)-k+V)v=0~~~\mathrm{and} ~~~\mathcal{P}' v = (P(i\nabla)+k-V^T)v=0. 
\]
\newdefi{Proposition} \label{weaksol}\emph{ If $w\in H^1(D)^8$ solves the vector Schr\"{o}dinger equation $ [-(\Delta +k^2) +Q]w=0$ weakly in a bounded Lipschitz domain $D$, i.e., $w$ satisfies
\begin{equation}
\int_{\mathbb{R}^3} \sum_{j=1}^8\big(\nabla w_j\cdot \nabla \varphi_j\big) - k^2 w\cdot \varphi \,dx + \langle Q w, \varphi \rangle = 0 \label{eq:schrw}
\end{equation}
for all $\varphi \in C_c^\infty(D)^8$, then $v=\mathcal{P}'w$ is a weak solution to 
\begin{equation}
\mathcal{P} v=0~~~\mathrm{in}~D, \label{eq:ellv}
\end{equation}
and $v\in H^1(D)^8$.}\vspace{2mm}

{\it Proof}. In order to show that $v$ is a weak solution to \eqref{eq:ellv}, we need to show that 
\begin{eqnarray*}
 0 &=&\int_{\mathbb{R}^3} v \cdot \mathcal{P}' \varphi\,dx=\int_{\mathbb{R}^3} \mathcal{P}' w \cdot \mathcal{P}'\varphi\,dx\\
 &=& \int_{\mathbb{R}^3} P(i\nabla)w \cdot P(i\nabla)\varphi + k^2 w\cdot \varphi + k w\cdot P(i\nabla)\varphi + P(i\nabla)w \cdot k \varphi \\
 &&\qquad - V^T w\cdot P(i\nabla)\varphi - P(i\nabla)w\cdot V^T \varphi - kw\cdot V^T\varphi - V^T w\cdot k \varphi + V^T w\cdot V^T \varphi\,dx
\end{eqnarray*}
for all $\varphi \in C_c^\infty(D)^8$, and we will do so by showing that this integral equals the left-hand side of \eqref{eq:schrw}. We first note that 
\[ \int_{\mathbb{R}^3} P(i\nabla)w\cdot P(i\nabla)\varphi\,dx = -\int_{\mathbb{R}^3} \sum_{l=1}^8 \nabla w_l \cdot \nabla \varphi_l\,dx. \]
The third and fourth terms cancel after an integration by parts, and it is straightforward to check that the last five terms give $-\langle Q w, \varphi\rangle$.

In order to see that $v\in H^1(D)^8$, we first note that certainly $v\in L^2(D)^8$; hence, by \eqref{eq:ellv}, we see that $P(i\nabla)v \in L^2(D)^8$, and by the ellipticity of $P(i\nabla)$, this implies $v\in H^1(D)^8$.\bewende
We also have the following analog when we switch the roles of $\mathcal{P}$ and $\mathcal{P}'$: \vspace{2mm}

\newdefi{Proposition} \label{weaksol2} \emph{If $w\in H^1(D)^8$ is a weak solution to $ [-(\Delta +k^2) +\tilde{Q}]w=0$ in $D$, i.e., $w$ satisfies
\[
\int_{\mathbb{R}^3} \sum_{j=1}^8\big(\nabla w_j\cdot \nabla \varphi_j\big) - k^2 w\cdot \varphi \,dx + \langle \tilde{Q}\, w, \varphi \rangle = 0
\]
for all $\varphi \in C_c^\infty(D)^8$, then $v=\mathcal{P}w$ is a weak solution to 
\[
\mathcal{P}'v=0~~~\mathrm{in}~D,
\]
and $v\in H^1(D)^8$.}\bewende
\subsection{Integral formula}\label{sec:intform}
Suppose now we have two sets of parameters, $ \mu_j, \varepsilon_j, \sigma_j \in C^{0,1}(\overline{\Omega})$, $j=1,2$, that have the same Cauchy data set on $\partial \Omega$, and such that on $\partial \Omega$, $\mu_1=\mu_2,\,\varepsilon_1=\varepsilon_2$, and $\sigma_1=\sigma_2$. Thus, we can perform a Whitney extension of the parameters to obtain Lipschitz continuous functions on the whole space such that $\mu_1=\mu_2,\,\varepsilon_1=\varepsilon_2$, and $\sigma_1=\sigma_2$ outside $\Omega$, and the Lipschitz constant on $\mathbb{R}^3$ depends only on the Lipschitz constant $c_o$ on $\overline{\Omega}$, cf. \cite[Section VI.2]{Stsingular}. We further require that outside of a sufficiently large ball $\Omega'$ containing $\Omega$, $\mu_1=\mu_2=\mu_o$, $\varepsilon_1=\varepsilon_2=\varepsilon_o$, and $\sigma_1=\sigma_2=0$. 

Set $\gamma_j = \varepsilon_j+i\sigma_j/ \omega$, and let $Q_j$ denote the weak potential with parameters $\mu_j,\,\gamma_j$, $j=1,2$. We also denote 
$$\mathcal{P}_j = P(i\nabla) -k + V(\mu_j,\gamma_j),~~~~\mathcal{P}'_j = P(i\nabla) +k - V^T(\mu_j,\gamma_j).$$ 
We start out by deriving the following integral formula that will be the starting point for the uniqueness proof. Note that for technical reasons, the functions we will construct will be solutions to the respective equations in the bigger bounded domain $\Omega' \supset \Omega$, but we have the Cauchy sets on $\partial \Omega$ available. \vspace{2mm}

\newdefi{Proposition}\label{intform} \emph{Let $w_1\in H^1_{loc}(\mathbb{R}^3)^8$ be a weak solution to the equation $ [-(\Delta +k^2) +Q_1]w_1=0$ in the bounded domain $\Omega'$, i.e., $w$ satisfies \eqref{eq:schrw} with $Q=Q_1$ for all $\varphi \in C_c^\infty(\Omega')^8$, and assume that $v_1= \mathcal{P}_1'w_1$ 
has vanishing first and last components. Furthermore, let $v_2\in H^1_{loc}(\mathbb{R}^3)^8$ satisfy 
\begin{equation}
\mathcal{P}_2' v_2 = 0~~~\mathrm{in}~\Omega'.  \label{eq:eqv2}
\end{equation}
Then, if the Cauchy data sets $C(\mu_1,\varepsilon_1,\sigma_1;\omega)$ and $C(\mu_2,\varepsilon_2,\sigma_2;\omega)$ on $\partial \Omega$ are equal, the following integral identity holds:
\begin{equation}
\big\langle (Q_2-Q_1) w_1, v_2 \big\rangle =0.\label{eq:intform}
\end{equation}}
{\it Proof.} As in the proof of Proposition \ref{weaksol}, we see that for $l=1,2$, 
\[
\int_{\mathbb{R}^3} \sum_{j=1}^8\big(\nabla w_{1,j}\cdot \nabla \varphi_j\big) - k^2 w_1\cdot \varphi \,dx + \langle Q_l w_1, \varphi \rangle = \int_{\mathbb{R}^3} \mathcal{P}_l' w_1 \cdot \mathcal{P}_l' \varphi \,dx
\]
for all $\varphi \in C^\infty_c(\mathbb{R}^3)^8$; by a density argument, the same holds for $\varphi \in H^1_{loc}(\mathbb{R}^3)^8$, so we may let $\varphi=v_2$. Subtracting the equation for $l=1$ from that for $l=2$, and using the fact that the two sets of parameters agree outside $\Omega$, we get 
\[ \big\langle (Q_2-Q_1) w_1, v_2 \big\rangle = \int_\Omega \mathcal{P}_2' w_1 \cdot \mathcal{P}'_2 v_2\, dx - \int_\Omega \mathcal{P}_1' w_1 \cdot \mathcal{P}'_1 v_2\, dx = -\int_\Omega v_1 \cdot \mathcal{P}'_1 v_2\, dx, \]
by the definition of $v_1$ and the assumption on $v_2$, by which the first integral vanishes. 
By construction, the assumption that $v_{1,1}=v_{1,4}=0$ guarantees that $u= (0,E,H,0):= (0, \gamma_1^{-1/2} v_{1,2}, \mu_1^{-1/2} v_{1,3}, 0)$ is a solution to Maxwell's equations in $\Omega'$ with parameters $\mu_1$ and $\gamma_1$,
\begin{equation}
 \nabla \wedge E - i\, \omega\, \mu_1\, H =0, \qquad \qquad \nabla \wedge H + i\, \omega\, \gamma_1\, E =0. \label{eq:me1}
\end{equation}
Thus, integrating by parts and using the fact that $\mathcal{P}_1v_1=0$ weakly, we obtain 
\begin{multline}
 -\int_\Omega v_1 \cdot \mathcal{P}'_1 v_2\, dx = \int_{\partial\Omega} (v_{1,2} \cdot \nu) v_{2,1} + v_{1,2}\cdot (\nu\wedge v_{2,3}) - v_{1,3}\cdot (\nu \wedge v_{2,2}) + (v_{1,3} \cdot \nu) v_{2,4} \,d S \\
 =  \int_{\partial\Omega} (\gamma_1 E \cdot \nu) \gamma_1^{-\frac{1}{2}} v_{2,1} - (\nu \wedge E) \cdot \gamma_1^{\frac{1}{2}} v_{2,3}+ (\nu\wedge H) \cdot \mu_1^{\frac{1}{2}} v_{2,2} + (\mu_1 H \cdot \nu) \mu_1^{-\frac{1}{2}} v_{2,4} \,d S. \label{eq:int1}
\end{multline}
Now we use the fact that $(\nu\wedge E|_{\partial \Omega},\nu\wedge H|_{\partial \Omega}) \in C_1=C_2$. This guarantees the existence of a solution $(\check{E},\check{H})\in \Hcurl(\Omega)^2$ to Maxwell's equations with parameters $\mu_2,\gamma_2$ in $\Omega$,
\begin{equation}
\nabla \wedge \check{E} - i\, \omega\, \mu_2\, \check{H} =0, \qquad \qquad \nabla \wedge \check{H} + i\, \omega\, \gamma_2\, \check{E} =0, \label{eq:me2}
\end{equation}
with the same Cauchy data. We set $g=(0, \gamma_2^{1/2} \check{E},\mu_2^{1/2}\check{H},0)$ almost everywhere in $\Omega$, and $g=v_1$ almost everywhere outside of $\Omega$. Then by \eqref{eq:eqv2}, we certainly have $\int_{\Omega} g\cdot \mathcal{P}_2' v_2 dx =0$, thus an analogous calculation to that above shows 
\begin{eqnarray}
 0 &=& \int_\Omega g \cdot \mathcal{P}'_2 v_2\, dx \notag \\
 &=&  -\int_{\partial\Omega} (\gamma_2 \check{E} \cdot \nu) \gamma_2^{-\frac{1}{2}} v_{2,1} - (\nu \wedge \check{E})\cdot \gamma_2^{\frac{1}{2}} v_{2,3} + (\nu\wedge\check{H} )\cdot \mu_2^{\frac{1}{2}} v_{2,2} + (\mu_2 \check{H} \cdot \nu) \mu_2^{-\frac{1}{2}} v_{2,4} \,d S. \label{eq:int2}
\end{eqnarray} 
Adding \eqref{eq:int2} to \eqref{eq:int1}, and using the fact that $\mu_1=\mu_2$ and $\gamma_1=\gamma_2$ on $\partial \Omega$, we obtain
\begin{eqnarray*}
\big\langle (Q_2-Q_1) w_1, v_2 \big\rangle &=&  \int_{\partial\Omega} (\gamma_1 E \cdot \nu) \gamma_1^{-\frac{1}{2}} v_{2,1} -(\gamma_2 \check{E} \cdot \nu) \gamma_2^{-\frac{1}{2}} v_{2,1} - (\mu_2 \check{H} \cdot \nu) \mu_2^{-\frac{1}{2}} v_{2,4}+ (\mu_1 H \cdot \nu) \mu_1^{-\frac{1}{2}}  v_{2,4} \,d S\\
 &=& \frac{1}{i\omega} \int_{\partial \Omega} \big\{[\nabla \wedge (\check{H}-H) ]\cdot \nu\big\} \gamma_1^{-\frac{1}{2}} v_{2,1}- \big\{[\nabla\wedge (\check{E}-E) ]\cdot \nu\big\} \mu_1^{-\frac{1}{2}} v_{2,4}\, d S,
 \end{eqnarray*}
where the last equality was obtained using Maxwell's equations \eqref{eq:me1} and \eqref{eq:me2}. Note that by our choice of $\check{H}$, $h:= \check{H}-H$ satisfies $\nu\wedge h = 0$ on $\partial \Omega$. Furthermore, for functions in $\Hcurl(\Omega')$, we have the identity $\nu\cdot \nabla\wedge h=-\mathrm{Div} \cdot (\nu \wedge h)$ in $H^{-1/2}(\partial \Omega)$, cf. \cite{MonkFEM}, so we find that $\nu\cdot (\nabla \wedge h) =0$, and similarly for $\check{E}-E$. This concludes the proof. \bewende
\subsection{A priori estimate}\label{sec:apriori}
We proceed by establishing an a priori estimate that is analogous to an estimate derived in \cite{CR2016}, and the proof will use some preliminary results from there. We first introduce the solution spaces, which were first introduced in \cite{HT2013} for scalar functions and are adapted to the structure of the equations we need to solve, as well as some auxiliary norms needed to prove the estimate: 
We denote by $p_\zeta$ the polynomial
\[ p_\zeta(\xi)= |\xi|^2 - 2 i \zeta\cdot \xi. \]
For $b\in \mathbb{R}$, we define $\dot{X}^{b}_{\zeta}$ to be the closure of the set of functions $w\in \mathcal{S}'(\mathbb{R}^3)^8$ for which 
\[ \|w\|_{\dot{X}^{b}_{\zeta}} = \bigg( \sum_{j=1}^8 \Big\| |p_\zeta |^b \hat{w}_j\Big\|_{L^2(\mathbb{R}^3)}^2\bigg)^{1/2} \]
is finite, and analogously $X^{b}_{\zeta}$ by the norm
\[ \|w\|_{X^{b}_{\zeta}} = \bigg( \sum_{j=1}^8 \Big\| \big(\big|\zeta\big|+\big|p_\zeta \big|\big)^b \hat{w}_j\Big\|_{L^2(\mathbb{R}^3)}^2\bigg)^{1/2}, \]
where $\hat{w}$ denotes the Fourier transform of $w$. 

For $M,\tau>0$, define the Fourier multiplier $m$ by
\[ m(\xi)= \left( M^{-1} \left| |\xi|^2-\tau^2\right|^2 + M^{-1} \tau^2 |\xi_n|^2 +M \tau^2\right)^{1/2}, \]
and for $u\in\mathcal{S}(\mathbb{R}^3)$, we define the norm
\[ \|u\|_{Y^b} = \| m^b\hat{u}\|_{L^2}. \]
For vector-valued functions we define the norm by 
\[ \|u\|_{Y^b}^2 = \sum_{k=1}^8{\| u_k\|_{Y^b}^2}. \]
The following estimate was proved in \cite[Lemma 2.3]{CR2016}: \vspace{2mm}

\newdefi{Lemma} \label{auxpriori} \emph{Let $\gamma$ be a Lipschitz continuous function that is constant outside a set of compact support, and let $A>1$ be such that 
\[ \| \gamma^{-1} \nabla\gamma\|_{L^\infty}<A. \]
Define $q=\gamma^{-1/2}\Delta \gamma^{1/2}$ in the weak sense, that is, for $\phi,\psi\in H^1_{loc}(\mathbb{R}^3)$
\[ \langle q\phi,\psi\rangle = \frac{1}{4} \int |\nabla \log \gamma|^2 \phi \psi - \frac{1}{2} \int \nabla\log \gamma \cdot \nabla(\phi\psi),\]
and furthermore, for a rotation $T$, let
\[ \langle T^*q \phi,\psi\rangle = \frac{1}{4} \int |\nabla \log \gamma(Tx)|^2 \phi(x) \psi(x)- \frac{1}{2} \int \nabla\log \gamma(Tx) \cdot \nabla( \phi(x) \psi(x)).\]
There is a constant $C$ such that for $M= C R^2 A^4$, and any $u\in \mathcal{S}(\mathbb{R}^3)$ with supp\,$u \subset \{ |x_3|<R\}$ and $\tau >8MR$,
\[ \|u\|_{Y^{1/2}} \lesssim \| (-\Delta + 2\tau \partial_{x_3} -\tau^2 +T^*q) u\|_{Y^{-1/2}}. \]
The implicit constant depends on $A$ and $R$. }\bewende
We use the following vector-valued generalization of this estimate with the weak matrix multiplier $Q$ of the form \eqref{eq:Qweak} (or $\tilde{Q}$ as in \eqref{eq:Qtilde}) instead of the scalar $q$: Recall that we set $\alpha=\nabla\log \gamma$ and $\beta=\nabla\log \mu$. Let $A>1$ be such that
\[ \max\big\{ \| \alpha \|_{L^\infty},\| \beta \|_{L^\infty}\big\} <A, \]
and let $Q$ be as in \eqref{eq:Qweak} (or \eqref{eq:Qtilde}). Then, as in the scalar case, there is a constant $C$ such that for $M= C R^2 A^4$, and any $u\in \mathcal{S}(\mathbb{R}^3)^8$ with supp\,$u \subset \{ |x_3|<R\}$ and $\tau >8MR$,
\[ \|u\|_{Y^{1/2}} \lesssim \| (-\Delta + 2\tau \partial_{x_3} -\tau^2 +T^*Q) u\|_{Y^{-1/2}}. \] 
We use this to establish the following estimate in the $X^b_{\zeta}$ norms, which is analogous to Proposition 2.4 in \cite{CR2016}. \vspace{2mm}

\newdefi{Proposition} \label{apriori} \emph{Let $\zeta = \real{\zeta} + i \imag \zeta \in \mathbb{C}^3$ such that $\real \zeta \perp \imag \zeta$ and $|\real \zeta|^2 = \tau^2= |\imag \zeta|^2 - k^2$. Furthermore, fix a constant $A> \max\{\|\alpha\|_{L^\infty},\|\beta\|_{L^\infty},1\}$. Then there exists an absolute constant $C$ such that for $\tau > C R^3 A^4$ and $|\zeta|$ sufficiently large,
\begin{equation}
\|u\|_{X^{1/2}_\zeta} \lesssim \| (-\Delta + 2\zeta \cdot \nabla +Q) \,u\|_{X^{-1/2}_\zeta}, \label{eq:apriori}
\end{equation}
provided that $u\in \mathcal{S}(\mathbb{R}^3)^8$ with supp\,$u \subset \{ |x|<R\}$. The implicit constant depends on $A$ and $R$.}\vspace{2mm}

{\it Proof.} Let $T$ be a rotation such that $\real \zeta = \tau T e_3$, where $e_3=(0,0,1)^T$. Let  $w\in \mathcal{S}(\mathbb{R}^3)^8$ with supp$\,w \subset \{ |x|<R\}$ and consider $v(x)= T^*w(x)=w(Tx)$. By Lemma \ref{auxpriori}, we have that for $M=CR^2A^4$ and $\tau>8MR$,
\begin{equation}
\|v\|_{Y^{1/2}} \lesssim \| (-\Delta + 2\tau \partial_{x_3} -\tau^2 +T^*Q) v\|_{Y^{-1/2}}. \label{eq:pfapriori}
\end{equation}
The right-hand function is 
\[  (-\Delta + 2\tau \partial_{x_3} -\tau^2 +T^*Q) w(Tx) = -T^*(\Delta w) + T^* 2 (\real \zeta \cdot \nabla) I_8 w - T^* |\real \zeta |^2 I_8 w + T^*Q w. \]
Writing $w(x) = e^{-i\imag \zeta \cdot x} u(x)$, with $u\in\mathcal{S}(\mathbb{R}^3)^8$ with supp\,$u\subset \{|x|\leq R\}$, one can easily check that
\[ T^*[(\Delta  + 2 (\real \zeta \cdot \nabla)   - |\real \zeta |^2  + Q )w] =  T^*[e^{-i\imag \zeta \cdot x} (-\Delta +k^2 +2\zeta\cdot \nabla +Q) u], \]
and thus \eqref{eq:pfapriori} yields 
\begin{equation}
\|T^*( e^{-i\imag \zeta \cdot x} u )\|_{Y^{1/2}} \lesssim \| T^*[e^{-i\imag \zeta \cdot x} (-\Delta +k^2 +2\zeta\cdot \nabla +Q) u]\|_{Y^{-1/2}}.\label{eq:estY}
\end{equation}
In order to compute these norms and relate them to the $X^b_{\zeta}$ norms, we need the Fourier transform of a function of the form $T^*( e^{-i\imag \zeta \cdot x} f )$, for $f\in X^{b}_\zeta$. A quick computation shows that 
\[ \mathcal{F}[T^*( e^{-i\imag \zeta \cdot x}f)] (\xi) = \mathcal{F}[e^{-i\imag \zeta \cdot x}f](T\xi) = \mathcal{F}[f](T\xi + \imag \zeta), \]
thus,
\[
\|T^*( e^{-i\imag \zeta \cdot x} f )\|_{Y^{b}} ^2 = \| m^{b} \mathcal{F}[T^*( e^{-i\imag \zeta \cdot x}f)] \|_{L^2}^2 = \int m(\xi)^{2b} | \hat{f}(T\xi + \imag \zeta)|^2 d\xi. \]
Now a change of variables and the definition of the multiplier $m$ show that the right hand integral is proportional to $ \| f\|_{X^b_\zeta}^2$, and using this on both sides of \eqref{eq:estY} gives the estimate
\begin{equation*}
\|u\|_{X_\zeta^{1/2}} \lesssim \|(-\Delta +k^2 +2\zeta\cdot \nabla +Q) u\|_{X_\zeta^{-1/2}}.
\end{equation*}
Applying the triangle inequality on the right-hand side and using the definition of the norms, we obtain 
\[
\|u\|_{X_\zeta^{1/2}} \lesssim \|(-\Delta+2\zeta\cdot \nabla +Q)\, u\|_{X_\zeta^{-1/2}} + k^2 |\zeta|^{-2} \| u \|_{X^{1/2}_\zeta}, \]
and the last term on the right can be absorbed into the left-hand side if $|\zeta|$ is large enough, which finishes the proof. \bewende
\subsection{Construction of CGO solutions}\label{sec:cgo}
We now proceed to construct solutions $w_1$ and $v_2$ to plug into the integral formula we derived in Section \ref{sec:intform}.
In what follows, we will consider the following localized spaces, in view of constructing solutions in the domain $\Omega'$:
\[ X^b_\zeta(\Omega')=\big\{ u|_{\Omega'}: u \in X^b_\zeta\big\},\quad b>0,  \]
with the norm 
\[ \|u\|_{X^b_{\zeta}(\Omega')} = \inf\big\{ \|v\|_{X^b_\zeta} : u=v|_{\Omega'} \big\}, \] 
as well as
\[ X^b_{\zeta,c}(\Omega')=\big\{ u \in X^b_\zeta: \mathrm{supp}\, u \subset \overline{\Omega'}\big\},\quad b \in \mathbb{R},  \] with the norm of $X^b_\zeta$. Note that $X^b_{\zeta,c}(\Omega')$ is a Hilbert space, and we can define $X^{-b}_{\zeta}(\Omega')$ to be its dual space.

\subsubsection{Construction of $w_1$}

We first choose a vector $\zeta_1 \in \mathbb{C}^3$ in the following way: Fix $\rho \in \mathbb{R}^3$, and choose unit vectors $\eta_1,\,\eta_2 \in \mathbb{R}^3$ such that $\{ \rho, \,\eta_1,\,\eta_2\}$ is an orthogonal basis of $\mathbb{R}^3$. Let $s\in \mathbb{R}$ with $s\geq 1$. Set
\[ \zeta_1 = -\sqrt{s^2+\frac{|\rho|^2}{4}} \,\eta_1 + i \left( \frac{1}{2}\rho - \sqrt{s^2+k^2} \,\eta_2\right). \]
Note that we have  
$\zeta_1\cdot \zeta_1=-k^2$. 
We now look for weak solutions of the form $w_1(x) = e^{\zeta_1\cdot x} (A_{\zeta_1} + R_{\zeta_1}(x))$ to the vector Schr\"odinger equation
\begin{equation}
 \big(-(\Delta + k^2) + Q_1\big)w_1=0, \label{eq:w_1}
 \end{equation}
where $A_{\zeta_1}$ is a constant 8-vector, and $R_{\zeta_1} \in X^{1/2}_{\zeta_1}$ compactly supported. 
In view of Proposition \ref{intform}, we need to guarantee that eventually $v=\mathcal{P}_1' w_1$ has vanishing first and last components. We want to facilitate this by choosing $A_{\zeta_1}$ such that the constant parts in the first and last components of $v$ vanish, and by imposing zero boundary conditions on $R_{\zeta_1}$ and extending the function by zero outside the domain so as to obtain a compactly supported function. This is the reason why we construct solutions in $\Omega'$, outside of which $\mu_1=\mu_2$ and $\gamma_1=\gamma_2$, so that in particular, $Q_1$ vanishes outside $\Omega'$. 

Plugging the ansatz for $w_1$ into \eqref{eq:w_1}, we obtain the following equation for $R_{\zeta_1}$, satisfied in the weak sense
\begin{equation}
 (-\Delta - (2\zeta_1 \cdot \nabla) I_8 + Q_1)R_{\zeta_1}=-Q_1 A_{\zeta_1}, \label{eq:R_1}
 \end{equation}
 that is, for any $\varphi  \in (\mathcal{S}(\mathbb{R}^3))^8$ with $\mathrm{supp}\,\varphi \subset \Omega'$, 
\begin{equation}
\left\langle R_{\zeta_1}, (-\Delta + (2 \zeta_1 \cdot \nabla)I_8 + Q_1) \varphi \right \rangle = -\left\langle Q_1 A_{\zeta_1}, \varphi\right\rangle. \label{eq:weak_R}
\end{equation}

We now want to find solutions to \eqref{eq:weak_R} that lie in the space $X^{1/2}_{\zeta_1,c}(\Omega')$, and we do so by defining a suitable linear functional on $X^{-1/2}_{\zeta_1}(\Omega') = (X^{1/2}_{\zeta_1,c}(\Omega'))^*$ which will be represented by $R_{\zeta_1}\in X^{1/2}_{\zeta_1,c}(\Omega')$.

\newdefi{Lemma}\label{R1} \emph{For given $\zeta_1$ as above, there exists a solution $R_{\zeta_1}\in X^{1/2}_{\zeta_1,c}(\Omega')$ to \eqref{eq:weak_R}.
Furthermore, $R_{\zeta_1}$ satisfies the averaged estimate 
\begin{equation}
\frac{1}{\lambda} \int_{S^1}\int_\lambda^{2\lambda} \|R_{\zeta_1}\|^2_{X^{1/2}_{\zeta_1}} \,ds\,d\eta_1 = o\big(\bm{1}(\lambda)\big),~~~\lambda\rightarrow \infty.
\label{eq:R1_avgd}
\end{equation}}
{\it Proof.}  We define the linear subspace
\[ \mathcal{L} = \Big\{ (-\Delta + (2 \zeta_1\cdot \nabla)I_8 + Q_1)\, \varphi: \varphi \in (\mathcal{S}(\mathbb{R}^3))^8,\,\mathrm{supp}\, \varphi \subset \Omega'\Big\} \subset X^{-1/2}_{\zeta_1}(\Omega'), \]
and the linear functional $L$ on $\mathcal{L}$ by
\[ Lv = -\left\langle Q_1 A_{\zeta_1}, \varphi\right\rangle,  \]
where $\varphi$ is such that $v=(-\Delta + (2 \zeta_1\cdot \nabla)I_8 + Q_1) \varphi$. Using the estimate \eqref{eq:apriori}, we get
\[ | Lv| \leq \| Q_1 A_{\zeta_1} \|_{X^{-1/2}_{\zeta_1}}   \| \varphi\|_{X^{1/2}_{\zeta_1}}   \lesssim \| Q_1 A_{\zeta_1} \|_{X^{-1/2}_{\zeta_1}} \| (-\Delta + (2\zeta_1\cdot \nabla)I_8 +Q_1) \varphi\|_{X^{-1/2}_{\zeta_1}}  = \| Q_1 A_{\zeta_1} \|_{X^{-1/2}_{\zeta_1}} \|v\|_{X^{-1/2}_{\zeta_1}} \]
for $v\in \mathcal{L}$, which shows that $L$ is a well-defined and bounded functional on $\mathcal{L}$. Thus, by the Hahn Banach Theorem, there exists an extension, still denoted by $L$, to $X^{-1/2}_{\zeta_1}(\Omega')$ that has the same operator norm.

Now $(X^{-1/2}_{\zeta_1}(\Omega'))^*\cong X^{1/2}_{\zeta_1,c}(\Omega')$ and therefore there is some $R_{\zeta_1}\in X^{1/2}_{\zeta_1,c}(\Omega')$ such that 
\[ Lv = \left\langle R_{\zeta_1}, v \right \rangle \]
for all $v\in X^{-1/2}_{\zeta_1}(\Omega')$. In particular, $R_{\zeta_1}$ satisfies \eqref{eq:weak_R}. We also have 
\begin{equation}
\| R_{\zeta_1}\|_{X^{1/2}_{\zeta_1}} \lesssim \| Q_1 A_{\zeta_1} \|_{X^{-1/2}_{\zeta_1}}. \label{eq:r-qa-est}
\end{equation}
From this we will obtain the averaged estimate for $R_{\zeta_1}$; first, we derive an estimate for $\| Q_1 A_{\zeta_1} \|_{X^{-1/2}_{\zeta_1}}$ using duality: let $\varphi \in X^{1/2}_{\zeta_1}$, and let $\chi \in C_o^{\infty}(\mathbb{R}^3)$ with $\chi(x)=1$ for $x\in \Omega'$. Since $Q_1$ is compactly supported on this set, we can estimate
\begin{equation}
\left| \left\langle Q_1 A_{\zeta_1} \,,\, \varphi \right\rangle \right| \lesssim \| \chi\varphi \|_{L^2} + \| FA_{\zeta_1}\|_{X^{-1/2}_{\zeta_1}} \|\varphi\|_{X^{1/2}_{\zeta_1}}, \label{eq:qa-norm}
\end{equation}
where the first term was obtained using boundedness of the material parameters and their first order derivatives, and (in the following, $\alpha= \nabla\log \gamma_1$, $\beta = \nabla \log \mu_1$, and $\kappa=\omega (\mu_1\gamma_1)^{1/2}$; $2(\nabla \beta)_S=\nabla \beta + (\nabla \beta)^T$)
\[
FA_{\zeta_1} = \left(\begin{array}{c}
(\nabla \cdot \alpha) A_{\zeta_1,1}\\
(2\big(\nabla \beta)_S - \nabla \cdot \beta\big) A_{\zeta_1,2}\\
(2\big(\nabla \alpha)_S - \nabla \cdot \alpha\big) A_{\zeta_1,3}\\
(\nabla \cdot \beta) A_{\zeta_1,4} \end{array} \right) \in X^{-1/2}_{\zeta_1}.
\]
Now, for the first term in \eqref{eq:qa-norm}, we can apply an adaptation of the estimate (5) from Lemma 2.2 in \cite{HT2013}, as well as the the estimate $\|\varphi\|_{\dot{X}^{1/2}_{\zeta_1}}\lesssim \|\varphi\|_{X^{1/2}_{\zeta_1}}$ to get
\[
\| \chi\varphi \|_{L^2} \lesssim s^{-1} \|\varphi\|_{\dot{X}^{1/2}_{\zeta_1}}\lesssim s^{-1} \|\varphi\|_{X^{1/2}_{\zeta_1}}, \]
whence from \eqref{eq:qa-norm} we get
\begin{equation}
 \|Q_1A_{\zeta_1}\|_{X^{-1/2}_{\zeta_1}}  \lesssim s^{-1}+\|FA_{\zeta_1}\|_{X^{-1/2}_{\zeta_1}}.\label{eq:qa-aux}
\end{equation}
In order to deal with the latter term, we further adapt Lemma 3.1 from \cite{HT2013} to our situation and obtain 
\[
\frac{1}{\lambda} \int_{S^1}\int_\lambda^{2\lambda} \|FA_{\zeta_1}\|^2_{\dot{X}^{-1/2}_{\zeta_1}} ds\,d\eta_1 = o\big(\bm{1}(\lambda)\big),~~~\lambda\rightarrow \infty.
\]
Noting that $\|FA_{\zeta_1}\|^2_{X^{-1/2}_{\zeta_1}}\lesssim\|FA_{\zeta_1}\|^2_{\dot{X}^{-1/2}_{\zeta_1}}$, \eqref{eq:qa-aux} now yields an averaged estimate for $ \|Q_1A_{\zeta_1}\|_{X^{-1/2}_{\zeta_1}}$,
\begin{equation}
\frac{1}{\lambda} \int_{S^1}\int_\lambda^{2\lambda} \|Q_1A_{\zeta_1}\|^2_{X^{-1/2}_{\zeta_1}} ds\,d\eta_1 \lesssim\frac{1}{\lambda} \int_{S^1}\int_\lambda^{2\lambda} \Big(s^{-2}+\|FA_{\zeta_1}\|^2_{X^{-1/2}_{\zeta_1}} \Big)ds\,d\eta_1 = o\big(\bm{1}(\lambda)\big),~~~\lambda\rightarrow \infty,\label{eq:avgdQ}
\end{equation}
and finally, by \eqref{eq:r-qa-est},
\[
\frac{1}{\lambda} \int_{S^1}\int_\lambda^{2\lambda} \|R_{\zeta_1}\|^2_{X^{1/2}_{\zeta_1}} ds\,d\eta_1 = o\big(\bm{1}(\lambda)\big),~~~\lambda\rightarrow \infty.
\]
\bewende
As we discussed above, we now extend $R_{\zeta_1}$ by zero into a slightly bigger bounded set $\Omega'' \supset \Omega'$ and note that since $Q_1=0$ outside $\Omega'$, this extension is in fact a compactly supported solution to the equation \eqref{eq:R_1} in $\Omega''$. We use this fact to show that the CGO solutions we just obtained actually yield solutions to Maxwell's equations. Recall that this is the case if the first and last components of $v= \mathcal{P}_1' w_1 $ vanish. \vspace{2mm}

\newdefi{Proposition}\label{w1} \emph{Let $\zeta_1$ and $A_{\zeta_1}$ be as above, and let $w_1(x) = e^{\zeta_1\cdot x} (A_{\zeta_1} + R_{\zeta_1}(x))$. Then $w_1\in H_{loc}^1(\mathbb{R}^3)^8$ and $w_1$ satisfies 
\begin{equation*}
\int_{\mathbb{R}^3} \sum_{j=1}^8\big(\nabla w_{1,j}\cdot \nabla \varphi_j\big) - k^2 w_1\cdot \varphi \,dx + \langle Q_1 w_1, \varphi \rangle = 0
\end{equation*}
for all $\varphi \in C_c^{\infty}(\Omega'')$. Furthermore, if $A_{\zeta_1}$ satisfies 
\[ i\zeta_1\cdot A_{\zeta_1,2}+kA_{\zeta_1,1} = i\zeta_1\cdot A_{\zeta_1,3}+kA_{\zeta_1,4} =0, \]
then if $|\zeta_1|$ is sufficiently large, $v=\mathcal{P}_1' w_1$ has vanishing first and last components in $\Omega''$. }\vspace{2mm}

{\it Proof}. Note first that $R_{\zeta_1} \in H^1(\mathbb{R}^3)^8$, since it is compactly supported in $\Omega''$. It follows that  $w_1\in H^1_{loc}(\mathbb{R}^3)^8$. 
Also, by construction, it is clear that $w_1$ solves \eqref{eq:w_1} in $\Omega''$. 

We proceed to prove the statement about $v=\mathcal{P}_1' w_1$. By Proposition \ref{weaksol}, the restriction of $v$ to $\Omega''$ belongs to $H^1(\Omega'')^8$ 
and $v$ satisfies $\mathcal{P}_1 v=0$ in $\Omega''$. Thus, we also have
\[
\mathcal{P}_1' \mathcal{P}_1 v =0
\]
weakly in $\Omega''$, and writing out the first and eighth component of this equation, we obtain weak decoupled equations (recall that we write the vector $v$ as $(v_1,v_2,v_3,v_4)$, where $v_1,v_4$ are scalars and $v_2,v_3$ are 3-vectors)
\begin{align}
-(\Delta + k^2) v_1 + (-\frac{1}{2} \nabla \cdot \beta + \frac{1}{4} |\beta|^2 - \theta) v_1 &= 0, \label{eq:vcomp1}\\
-(\Delta + k^2) v_4 + (-\frac{1}{2} \nabla \cdot \alpha + \frac{1}{4} |\alpha|^2 - \theta) v_4 &= 0.\notag
\end{align}
Our next goal is to show that $v_1$ and $v_4$ vanish if only we pick $A_{\zeta_1}$ suitably. Using the shape of $w_1$, we get
\begin{align*} 
v &= (P(i\nabla)+k-V(\mu_1,\gamma_1)^T) e^{\zeta_1\cdot x}(R_{\zeta_1}+A_{\zeta_1}) \\
&=e^{\zeta_1 \cdot x} \Big( \underbrace{(P(i\zeta_1)+k)A_{\zeta_1}}_{=const} + \big(iP(\nabla + \zeta_1)+k\big)R_{\zeta_1} -V^T(R_{\zeta_1}+A_{\zeta_1}) \Big).
\end{align*}
The equation for the first component $v_1$ reads
\begin{equation} 
v_1 =e^{\zeta_1 \cdot x} \Big( \underbrace{i\zeta_1\cdot A_{\zeta_1,2}+kA_{\zeta_1,1}}_{=const} + \underbrace{i(\nabla + \zeta_1)\cdot R_{\zeta_1,2} +k R_{\zeta_1,1} -(k-\kappa)(R_{\zeta_1,1}+A_{\zeta_1,1}) - \frac{i}{2} \beta\cdot (R_{\zeta_1,2}+A_{\zeta_1,2}) }_{=S}  \Big). \label{eq:v1}
\end{equation}
Furthermore, the last row is
\begin{equation*} 
v_4 =e^{\zeta_1 \cdot x} \Big( \underbrace{i\zeta_1\cdot A_{\zeta_1,3}+kA_{\zeta_1,4}}_{=const} + i(\nabla + \zeta_1)\cdot R_{\zeta_1,3} +k R_{\zeta_1,4} -(k-\kappa)(R_{\zeta_1,4}+A_{\zeta_1,4}) - \frac{i}{2} \alpha\cdot (R_{\zeta_1,3}+A_{\zeta_1,3})  \Big).
\end{equation*}
By the assumption on $A_{\zeta_1}$, the constant parts in both equations vanish. Thus, $v_1$ is of the form $v_1= e^{\zeta_1\cdot x} S$, with $S$ the non-constant part in \eqref{eq:v1}. Note that on $\partial \Omega''$ we have $k=\kappa$, $\beta =\nabla \log \mu_1 =0$, and $R_{\zeta_1}$ is compactly supported in $\Omega''$, therefore $v_1=0$ on $\partial \Omega''$. So $v_1\in H^1_0(\Omega'')$, 
and it follows that $v_1=0$ in $\Omega''$ by uniqueness of the solution to the Schr\"{o}dinger equation \eqref{eq:vcomp1} in $H^1_0(\Omega'')$. The fact that $v_4=0$ is proved analogously.\bewende
\newdefi{Remark} We will use the following choice of $A_{\zeta_1}$ that satisfies the condition of Proposition \ref{w1}: 
\[
A_{\zeta_1} = \frac{\sqrt{2}}{|\zeta_1|} \big(\zeta_1\cdot a, i k\,a, i k\,b, \zeta_1\cdot b \big)^T, 
\]
where $a$ and $b$ are constant vectors  in $\mathbb{R}^3$. Below, we will be using two different choices for these vectors, letting either $a=\eta_1$ and $b=0$, or $a=0$ and $b=|\rho|^{-1} \eta_2\wedge\rho$. 

\subsubsection{Construction of $v_2$}
We now construct the function $v_2$ that weakly solves 
\[
\mathcal{P}_2' v_2 = (P(i\nabla) + k - V(\mu_2,\gamma_2)^T)v_2=0~~~\mathrm{in}~ \Omega'.
\]
We first pick $\zeta_2\in \mathbb{C}^3$ in the following way: We take the orthogonal basis $\{ \rho, \,\eta_1,\,\eta_2\}$ and $s\in \mathbb{R}$ with $s\geq 1$ as above for $\zeta_1$, and set
\[ \zeta_2 = \sqrt{s^2+\frac{|\rho|^2}{4}} \eta_1 + i \left( \frac{1}{2}\rho + \sqrt{s^2+k^2} \eta_2 \right). \]
Note that this choice was made such that $\zeta_2\cdot \zeta_2 = -k^2$ and $\zeta_1+\zeta_2=i \rho$; the latter will be exploited when we plug the CGO solutions into the integral formula. We now want to find a solution of the form
\[ v_2(x) = e^{\zeta_2\cdot x}\big(B_{\zeta_2}+ S_{\zeta_2}(x)\big), \]  
with $B_{\zeta_2}$ a constant vector in $\mathbb{C}^8$ and $S_{\zeta_2}\in X^{1/2}_{\zeta_2}(\Omega')$, small on average. Analogously to the construction of $v$ from $w_1$, we start out by first finding a weak solution to a vector Schr\"{o}dinger equation.\vspace{2mm}

\newdefi{Lemma} \label{w2} \emph{Let $\zeta_2$ be as above, and let $A_{\zeta_2}$ be a constant vector with $A_{\zeta_2,1}=A_{\zeta_2,4}=0$. Then for $|\zeta_2|$ sufficiently large, there is a solution $w_2\in H^1(\Omega')^8$ of the form $w_2(x)= e^{\zeta_2\cdot x} (A_{\zeta_2} + R_{\zeta_2}(x))$ to 
\begin{equation}
[-(\Delta +k^2)+\tilde{Q}_2 ]w_2=0~~~\mathrm{in}~\Omega',\label{eq:w2}
\end{equation}
with $R_{\zeta_2}\in X^{1/2}_{\zeta_2,c}(\Omega')$, such that $w_{2,1}=w_{2,4}=0$. }\vspace{2mm}

{\it Proof.} Plugging the ansatz for $w_2$ into \eqref{eq:w2}, we find that $R_{\zeta_2}$ needs to weakly satisfy the equation
\[
\big(-\Delta - 2\zeta_2 \cdot \nabla + \tilde{Q}_2\big)R_{\zeta_2}=-\tilde{Q}_2 A_{\zeta_2}
\]
in $\Omega'$. The existence of such $R_{\zeta_2} \in X^{1/2}_{\zeta_2,c}(\Omega')$ follows by the same argument as in proving Lemma \ref{R1}. Furthermore, the structure \eqref{eq:Qtilde} of $\tilde{Q}_2$ shows that the first and last components of the equation decouple, so that the fact that $A_{\zeta_2,1}=A_{\zeta_2,4}=0$ yields $R_{\zeta_2,1}=R_{\zeta_2,4}=0$. Therefore, we get $w_{2,1}=w_{2,4}=0$. The fact that $w_2\in H^1(\Omega')^8$ follows since $R_{\zeta_2} \in H^1(\Omega')^8$. 
\bewende
In the following, we will use these specific choices of $A_{\zeta_2}$, adapted to the choices made for $A_{\zeta_1}$:
Let $a$ and $b$ be unit vectors chosen above, that is, either $a=\eta_1$ and $b=0$, or $a=0$ and $b=|\rho|^{-1} \eta_2\wedge\rho$. For either of these, we let $A_{\zeta_2}=-\frac{\sqrt{2}}{|\zeta_2|} (0,a,b,0)^T$. \vspace{2mm} 

\newdefi{Proposition} \label{v2} \emph{Let $\zeta_2$  and $A_{\zeta_2}$ be as above and let $w_2$ be the function constructed in Lemma \ref{w2}. Then $v_2=\mathcal{P}_2 w_2$ belongs to $H^1(\Omega')^8$ and is a weak solution to 
\begin{equation}
\mathcal{P}_2' v_2 =0 ~~in~\Omega'.\label{eq:v2eq}
\end{equation}
We can furthermore extend $v_2$ to $H^1_{loc}(\mathbb{R}^3)^8$ and write $v_2(x)=e^{\zeta_2\cdot x}(B_{\zeta_2}+ S_{\zeta_2}(x))$, with constant $B_{\zeta_2}$ and $S_{\zeta_2}\in X^{1/2}_{\zeta_2}$, where $S_{\zeta_2}$ satisfies }
\begin{equation}
\frac{1}{\lambda} \int_{S^1}\int_\lambda^{2\lambda} \|S_{\zeta_2}\|^2_{X^{1/2}_{\zeta_2}} ds\,d\eta_1 = o(1),~~~\lambda\rightarrow \infty.
\label{eq:S2_avgd}
\end{equation}

{\it Proof.} The fact that $v_2\in H^1(\Omega)^8$ and $v_2$ solves \eqref{eq:v2eq} are the statements of Proposition \ref{weaksol2}. 
Writing $v_2(x)=e^{\zeta_2\cdot x}(B_{\zeta_2}+ S_{\zeta_2}(x))$, we have 
\begin{gather*} 
B_{\zeta_2}=(P(i\zeta_2)-k)A_{\zeta_2} = -\frac{\sqrt{2}}{|\zeta_2|} \big(i\zeta_2\cdot a, \,i\zeta_2 \wedge b-k\,a, -i\zeta_2\wedge a - k\,b,\, i \zeta_2\cdot b \big)^T,\\
S_{\zeta_2} = i P(\nabla + \zeta_2) R_{\zeta_2} - k R_{\zeta_2} + V(\mu_2,\gamma_2)(A_{\zeta_2} + R_{\zeta_2}).
\end{gather*}
Applying $\mathcal{P}_2$ 
to \eqref{eq:v2eq}, we find that $S_{\zeta_2}$ weakly satisfies the equation
\begin{equation}
 (-\Delta -2\zeta_2\cdot \nabla +Q_2) S_{\zeta_2} = -Q_2 B_{\zeta_2} ~~\mathrm{in}~\Omega', \label{eq:de-s2}
 \end{equation}
where $-Q_2 B_{\zeta_2}\in X^{-1/2}_{\zeta_2}$. Using the fact that $R_{\zeta_2}=0$ on $\partial \Omega'$, and that $Q_2$ and $V(\mu_2,\gamma_2)$ vanish outside $\Omega'$, we can again extend $S_{\zeta_2}=0$ to the slightly bigger domain $\Omega''$, and see that this extension solves \eqref{eq:de-s2} in $\Omega''$ with zero boundary condition. This elliptic equation has a unique solution in $H^1_0(\Omega'')^8$, and since the corresponding extension of $v_2$ to $\Omega''$ belongs to $H^1(\Omega'')^8$, we find that $S_{\zeta_2}\in H^1_0(\Omega'')^8$, and hence $S_{\zeta_2}$ is this unique solution.

On the other hand, we can employ the method used to find $R_{\zeta_2}$ to solve \eqref{eq:de-s2} in $X^{1/2}_{\zeta_2,c}(\Omega'')$ and obtain a solution $\tilde{S} \in X^{1/2}_{\zeta_2,c}(\Omega'') \subset H^1_0(\Omega'')^8$ that satisfies the estimate
\[ \| \tilde{S}\|_{X^{1/2}_{\zeta_2}} \lesssim \| Q_2 B_{\zeta_2} \|_{X^{-1/2}_{\zeta_2}}. \]
But uniqueness in  $H^1_0(\Omega'')^8$ now shows that $\tilde{S}=S_{\zeta_2}$, and the above estimate can be used to prove \eqref{eq:S2_avgd} in the same manner as the corresponding estimate for $R_{\zeta_1}$ in Lemma \ref{R1} (note that $| B_{\zeta_2}| = O(1)$ as $|\zeta_2|\rightarrow \infty$). The extension to $H^1_{loc}(\mathbb{R}^3)^8$ is performed by setting $S_{\zeta_2}=0$ outside $\Omega''$. \bewende
\subsection{Uniqueness of the parameters} \label{sec:unique}
We now plug the solutions $w_1$ and $v_2$ constructed above into the integral formula \eqref{eq:intform} to get
\[ \big\langle (Q_2-Q_1)e^{\zeta_1\cdot x}(R_{\zeta_1}+A_{\zeta_1}) , e^{\zeta_2\cdot x}(S_{\zeta_2}+B_{\zeta_2}) \big\rangle =0. \]
By our choice of $\zeta_1$ and $\zeta_2$, $e^{\zeta_1\cdot x} e^{\zeta_2\cdot x}=e^{i \rho \cdot x}$. Recall also our choices
\begin{gather*}
A_{\zeta_1} = \frac{\sqrt{2}}{|\zeta_1|} \big(\zeta_1\cdot a, i k\,a, i k\,b, \zeta_1\cdot b \big)^T, \quad 
B_{\zeta_2} = -\frac{\sqrt{2}}{|\zeta_2|} \big(i\zeta_2\cdot a, \,i\zeta_2 \wedge b-k\,a, -i\zeta_2\wedge a - k\,b,\, i \zeta_2\cdot b \big)^T.
\end{gather*}
We define $A_1$ and $B_2$ as the limits of $A_{\zeta_1}$ and $B_{\zeta_2}$, respectively, as $s\rightarrow \infty$. That is,
\begin{align*}
A_{1} &=  \big(-(\eta_1+i\eta_2)\cdot a, \,0,\,0, -(\eta_1+i \eta_2)\cdot b \big)^T,\\
B_{2} &=  -i \big((\eta_1+i\eta_2)\cdot a, (\eta_1+i\eta_2)\wedge b,-(\eta_1+i\eta_2)\wedge a, (\eta_1+i \eta_2)\cdot b \big)^T.
\end{align*}
Explicitly, for $a=\eta_1$ and $b=0$ we get
\begin{gather*}
A_{1} =  \big(-1, \,0,\,0,\,0 \big)^T,\quad
B_{2} = \big(-i,\,0, \eta_1\wedge \eta_2 ,\,0),
\end{gather*}
and for $a=0$ and $b=|\rho|^{-1} \eta_2\wedge\rho$,
\begin{gather*} A_{1} = \big(0, \,0,\,0, -1 \big)^T,\quad B_{2} = \big(0,-\frac{\rho}{|\rho|},\,0,\,-i\big). \end{gather*}
Note that for any choice of $\rho$, the convergence rate for large $s$ is $|A_{\zeta_1} -A_1|+|B_{\zeta_2}-B_2| =\mathcal{O}(s^{-1})$.

Subtracting the term $ \langle (Q_2-Q_1)e^{i \rho\cdot x}A_1|B_2\rangle$ as well as adding and subtracting $ \langle (Q_2-Q_1)(R_{\zeta_1}+A_{\zeta_1}),e^{i \rho\cdot x}B_2\rangle$, we obtain the equation
\begin{multline} 
\qquad -\big \langle (Q_2-Q_1) e^{i \rho\cdot x} A_1,B_2\big\rangle =  \big\langle (Q_2-Q_1)(R_{\zeta_1}+A_{\zeta_1}),e^{i \rho\cdot x}(B_{\zeta_2}-B_2 +S_{\zeta_2})\big \rangle\\
 + \big \langle (Q_2-Q_1)B_2,e^{i \rho\cdot x}(R_{\zeta_1}+A_{\zeta_1}-A_1)\big\rangle,\qquad \label{eq:deltaQ}
\end{multline}
where we have used the symmetry of the operators $Q_i$. Our goal now is to show that the right-hand side of this equation tends to zero on average, as $s$ becomes large (note that the left-hand side does not depend on $s$), so that we get
\[ \big\langle (Q_2-Q_1)A_1 e^{i\rho\cdot x},B_2\big\rangle=0. \]
From this equation, we will be able to extract a set of two differential equations for the material parameters, using the two different choices for $A_1$ and $B_2$ specified above.\vspace{3mm}\\
We have already seen that for each $\rho$, $Q_i$ are bounded from $X_{\zeta_i}^{1/2}$ to $X_{\zeta_i}^{-1/2}$. Furthermore, $Q_1-Q_2$ is bounded from $X_{\zeta_1}^{1/2}$ to $X_{\zeta_2}^{-1/2}$: In order to see this, we need to show that for any $w\in X_{\zeta_1}^{1/2}$ and $v\in X_{\zeta_2}^{1/2}$,  
\[ \big|\big\langle (Q_2-Q_1)w,v\big\rangle\big| \lesssim \|w\|_{X_{\zeta_1}^{1/2}} \|v\|_{X_{\zeta_2}^{1/2}}. \]
But this follows from the definiton of $Q_2-Q_1$ and the argument used in the proof of Lemma 2.3 in \cite{HT2013}, since $\langle (Q_2-Q_1)w|v\rangle$ involves expressions of the form
\[ \int{q(x) w_i(x) v_k(x)\, dx}\qquad \mathrm{and} \qquad \int{\tilde{q}(x) \partial_{x_j}(w_i(x) v_k(x))dx}, \]
with the functions $q,\,\tilde{q}$ being Lipschitz or bounded functions. Note that \cite[Lemma 2.3]{HT2013} gives an estimate involving the $\dot{X}^b_{\zeta_i}$ norms; the estimate in the $X^b_{\zeta_i}$ norms follows using that for compactly supported functions we have
\[ \|f \|_{X^{-1/2}_{\zeta}} \lesssim  \|f \|_{\dot{X}^{-1/2}_{\zeta}},~~~\mathrm{and}~~~  \|f \|_{\dot{X}^{1/2}_{\zeta}} \lesssim  \|f \|_{X^{1/2}_{\zeta}}. \ \]
Letting $\chi \in C^\infty_0(\mathbb{R}^3)$ such that $\chi = 1$ on $\Omega'$, we can estimate \eqref{eq:deltaQ} by
\begin{multline}
\big|\big\langle (Q_2-Q_1)A_1 e^{i\rho\cdot x},\,B_2\big\rangle\big| \lesssim \big(\|(Q_2-Q_1) A_{\zeta_1} \|_{X^{-1/2}_{\zeta_2}} + \|R_{\zeta_1}\|_{X^{1/2}_{\zeta_1}}\big)\big( \|(B_{\zeta_2}-B_2) \chi\|_{X^{1/2}_{\zeta_2}}+\|S_{\zeta_2}\|_{X^{1/2}_{\zeta_2}}\big) \\
+ \|(Q_2-Q_1) B_{\zeta_2}\|_{X^{-1/2}_{\zeta_1}} \big( \|R_{\zeta_1}\|_{X^{1/2}_{\zeta_1}} + \|(A_{\zeta_1}-A_1) \chi\|_{X^{1/2}_{\zeta_1}}\big). \label{eq:absdeltaQ}
\end{multline}
Using Lemma 2.2 from \cite{HT2013}, we obtain 
\[ \left(\frac{1}{\lambda} \int_{S^1}\int_\lambda^{2\lambda} \|\chi (A_{\zeta_1}-A_1)\|_{X^{1/2}_{\zeta_1}}^2 ds\, d\eta_1\right)^{1/2} = O(\bm{1}(\lambda)),~~~\lambda \rightarrow \infty, \]
\[ \left( \frac{1}{\lambda} \int_{S^1}\int_\lambda^{2\lambda} \|\chi (B_{\zeta_2}-B_2)\|_{X^{1/2}_{\zeta_2}}^2 ds \,d\eta_1 \right)^{1/2}= O(\bm{1}(\lambda)),~~~\lambda \rightarrow \infty, \]
and thus if we take the average over $(s,\eta_1)\in [\lambda,2\lambda]\times S^1$ of \eqref{eq:absdeltaQ}, and use the Cauchy-Schwarz inequality as well as the estimates \eqref{eq:R1_avgd} and \eqref{eq:S2_avgd}, we obtain
\begin{multline}
\big|\big\langle (Q_2-Q_1)A_1 e^{i\rho\cdot x},\,B_2\big\rangle\big| \lesssim \Big( O\big(\bm{1}(\lambda)\big) + o\big(\bm{1}(\lambda)\big)\Big) \left( \Big( \frac{1}{\lambda} \int_{S^1}\int_{\lambda}^{2\lambda}  \|(Q_2-Q_1) A_{\zeta_1}\|_{X^{-1/2}_{\zeta_2}}^2 ds \,d\eta_1 \Big)^{1/2} + o\big(\bm{1}(\lambda)\big)\right)\\
+ \Big( O\big(\bm{1}(\lambda)\big) + o\big(\bm{1}(\lambda)\big)\Big) \left( \frac{1}{\lambda} \int_{S^1}\int_{\lambda}^{2\lambda}  \|(Q_2-Q_1) B_{\zeta_2}\|_{X^{-1/2}_{\zeta_1}}^2 ds \,d\eta_1 \right)^{1/2} \label{eq:avgdeltaQ}
\end{multline}
The same argument as in showing \eqref{eq:avgdQ} shows that the two integrals tend to zero as $\lambda \rightarrow \infty$, and since the left-hand side of \eqref{eq:avgdeltaQ} is independent of $\lambda$, we arrive at 
\[\big\langle (Q_2-Q_1)A_1 e^{i\rho\cdot x},\,B_2\big\rangle =0. \]
By plugging in each of our choices for $A_1$ and $B_2$, we obtain two differential equations for the coefficients. We let $\alpha_i=\nabla \log \gamma_i$, $\beta_i=\nabla \log \mu_i$, and $\theta_i=\omega^2\gamma_i\mu_i$ for $i=1,2$. Then, for $a=\eta_1$ and $b=0$, we get 
\[ 0 = \langle (Q_2-Q_1)A_1 e^{i\rho\cdot x},\,B_2\rangle 
=  \int{\frac{i}{4} \big(\alpha_2\cdot \alpha_2-4\theta_2 -(\alpha_1\cdot \alpha_1-4\theta_1)\big)e^{i\rho\cdot x}
-\frac{i}{2}(\alpha_2-\alpha_1)\cdot \nabla (e^{i \rho x}) dx}. \]
Since $\rho$ is an arbitrary vector, we conclude that the integrand must equal zero, and this gives the weakly satisfied equation, using $\alpha_2\cdot \alpha_2-\alpha_1\cdot\alpha_1 = (\alpha_2-\alpha_1)\cdot(\alpha_2+\alpha_1)$,
\[  (\alpha_2-\alpha_1)\cdot(\alpha_2+\alpha_1) -4(\theta_2-\theta_1) + 2 \nabla\cdot(\alpha_2-\alpha_1)=0. \]
Plugging in $a=0$ and $b= |\rho|^{-1} \eta_2\wedge\rho$ yields
\[ 0 = \int{ \frac{i}{4}\big(\beta_2\cdot \beta_2-4\theta_2 -(\beta_1\cdot \beta_1-4\theta_1)\big)e^{i\rho\cdot x}
-\frac{i}{2} (\beta_2-\beta_1)\cdot \nabla( e^{i \rho x})dx}, \]
hence
\[  (\beta_2-\beta_1)\cdot(\beta_2+\beta_1) -4(\theta_2-\theta_1) + 2 \nabla\cdot(\beta_2-\beta_1)=0. \]
In terms of $\gamma_i$ and $\mu_i$, we obtain the following system of equations.
\begin{align*}
\frac{1}{4} \left(\frac{\nabla \gamma_2}{\gamma_2} - \frac{\nabla \gamma_1}{\gamma_1} \right)\cdot \left(\frac{\nabla \gamma_2}{\gamma_2} + \frac{\nabla \gamma_1}{\gamma_1} \right)- \omega^2\big(\gamma_2\mu_2-\gamma_1\mu_1\big) +\frac{1}{2} \nabla\cdot \left(\frac{\nabla \gamma_2}{\gamma_2} - \frac{\nabla \gamma_1}{\gamma_1} \right) &= 0,\\
\frac{1}{4} \left(\frac{\nabla \mu_2}{\mu_2} - \frac{\nabla \mu_1}{\mu_1} \right)\cdot \left(\frac{\nabla \mu_2}{\mu_2} + \frac{\nabla \mu_1}{\mu_1} \right)- \omega^2\big(\gamma_2\mu_2-\gamma_1\mu_1\big) +\frac{1}{2} \nabla\cdot \left(\frac{\nabla \mu_2}{\mu_2} - \frac{\nabla \mu_1}{\mu_1} \right) &= 0.
\end{align*}
Some basic computations show that the above is equivalent to the following system,
\begin{align*}
-\Delta(\gamma_2^{1/2}-\gamma_1^{1/2}) + V(\gamma_2^{1/2}-\gamma_1^{1/2}) + a (\gamma_2^{1/2}-\gamma_1^{1/2}) + b (\mu_2^{1/2}-\mu_1^{1/2}) &= 0,\\
-\Delta(\mu_2^{1/2}-\mu_1^{1/2}) + W(\mu_2^{1/2}-\mu_1^{1/2}) + c (\mu_2^{1/2}-\mu_1^{1/2}) + d (\gamma_2^{1/2}-\gamma_1^{1/2}) &= 0.
\end{align*}
with
\[
V=\frac{\Delta(\gamma_1^{1/2}+\gamma_2^{1/2})}{\gamma_1^{1/2}+\gamma_2^{1/2}}, \qquad a=\bm{1}_{\Omega} \omega^2 \gamma_1^{1/2}\gamma_2^{1/2} (\mu_1+\mu_2), \qquad b=\bm{1}_{\Omega} \omega^2 \gamma_1^{1/2}\gamma_2^{1/2} (\gamma_1+\gamma_2) \frac{\mu_1^{1/2}+\mu_2^{1/2}}{\gamma_1^{1/2}+\gamma_2^{1/2}},
\]
\[
W=\frac{\Delta(\mu_1^{1/2}+\mu_2^{1/2})}{\mu_1^{1/2}+\mu_2^{1/2}}, \qquad c=\bm{1}_{\Omega}\omega^2 \mu_1^{1/2}\mu_2^{1/2} (\gamma_1+\gamma_2), \qquad d=\bm{1}_{\Omega}\omega^2 \mu_1^{1/2}\mu_2^{1/2}(\mu_1+\mu_2)\frac{\gamma_1^{1/2}+\gamma_2^{1/2}}{\mu_1^{1/2}+\mu_2^{1/2}}.
\]
Note that $\gamma_2^{1/2}-\gamma_1^{1/2}$ and $\mu_2^{1/2}-\mu_1^{1/2}$ are compactly supported functions in $H^1(\mathbb{R}^3)$. From the next result it will follow that both of these functions vanish and hence $\gamma_1=\gamma_2$ and $\mu_1=\mu_2$, which finishes the proof of Theorem \ref{thm1}.\vspace{2mm}

\newdefi{Lemma} \label{uniquecont} \emph{Suppose $f$ and $g$ are compactly supported functions in $H^1(\mathbb{R}^3)$. Then, if $f$ and $g$ satisfy the system of equations
\begin{align*}
-\Delta f + V f + af + bg&= 0,\\
-\Delta g + Wg +cg + df &= 0,
\end{align*}
the two functions vanish identically.}\vspace{2mm}

{\it Proof.} Let $\zeta \in \mathbb{C}^3$ satisfy $\zeta \cdot \zeta =0$, and set $u(x) = e^{-\zeta\cdot x} f(x)$ and $v(x)=e^{-\zeta\cdot x}g(x)$. Then $u$ and $v$ are compactly supported functions in $H^1(\mathbb{R}^3)$, and consequently also belong to $X^{1/2}_\zeta$.
Furthermore, $u$ and $v$ satisfy the equations
\begin{align}
-\Delta u +2 \zeta\cdot \nabla u + V u + au + bv&= 0, \label{eq:uv1}\\
-\Delta v + 2 \zeta\cdot \nabla v + Wv +cv + du &= 0.\label{eq:uv2}
\end{align}
Define the weak potential $Q$ by
\begin{align*}
\langle Q w, \varphi \rangle =& \int -\nabla\big( \gamma_1^{1/2}+\gamma_2^{1/2} \big) \cdot \nabla\left( \frac{w_1\varphi_1}{\gamma_1^{1/2}+\gamma_2^{1/2}} \right) + \big( a w_1 + b w_4 \big) \varphi_1 \notag \\
& \qquad - \nabla \big(\mu_1^{1/2}+\mu_2^{1/2}\big)\cdot \nabla \left( \frac{w_4\varphi_4}{\mu_1^{1/2}+\mu_2^{1/2}} \right) + \big(c w_4+ d w_1\big) \varphi_4 \,dx
\end{align*}
for $w=(w_1,0,0,w_4),\,\varphi=(\varphi_1,0,0,\varphi_4)$ in $H^1(\mathbb{R}^3)^8$, and note that $Q$ is compactly supported and bounded, and satisfies the conditions for the a priori estimate \eqref{eq:apriori}. 
Setting $w=(u,0,0,v)$, the system \eqref{eq:uv1}-\eqref{eq:uv2} can be written as
\[
(-\Delta+2 \zeta\cdot \nabla + Q) w  =0.
\]
A slight modification of the a priori estimate \eqref{eq:apriori} now shows that $w=0$, provided $|\zeta|$ is chosen large enough, which finishes the proof. \bewende

\begin{center}
{\sc Acknowledgment}
\end{center}
The author would like to thank her advisor, Ting Zhou, for her patient guidance and many helpful discussions.

\bibliographystyle{plain}
\bibliography{inverse.bib}
\vspace{5mm}

{\sc Department of Mathematics, Northeastern University. Boston, MA 02115. USA.} 

E-mail address: \href{mailto:pichler.mo@husky.neu.edu}pichler.mo@husky.neu.edu

\end{document}